\documentclass{birkjour}

\newcommand{\C}{{\mathbb C}}
\newcommand{\D}{{\mathbb D}}

\newcommand{\T}{{\mathbb T}}

\newcommand{\bmo}{{\rm BMO}}
\newcommand{\bmoa}{{\rm BMOA}}

\newcommand{\vmo}{{\rm VMO}}

\newcommand{\eps}{\varepsilon}

\newcommand{\f}{\frac}
\newcommand{\ov}{\overline}

\newcommand{\de}{\delta}

\newcommand{\ze}{\zeta}

\newcommand{\ph}{\varphi}

\title[Strong Reproducing Kernel Thesis and the Garsia Norm]
{The Strong Reproducing Kernel Thesis\\ 
for Hankel Operators and the Garsia Norm}
\author[K. M. Dyakonov]{Konstantin M. Dyakonov}

\address
{Departament de Matem\`atiques i Inform\`atica\\ 
Universitat de Barcelona, IMUB and BGSMath\\
Gran Via, 585\\ 
E-08007 Barcelona\\ 
Spain}

\address{\it \,\,\\
and}

\address{\quad\\
Instituci\'o Catalana de Recerca i Estudis Avan\c{c}ats (ICREA)\\ 
Pg. Llu\'is Companys, 23\\ 
E-08010 Barcelona\\ 
Spain}

\email{konstantin.dyakonov@icrea.cat}
\keywords{Hardy space, Hankel operator, reproducing kernel, bounded mean oscillation, Garsia norm, extreme point}
\subjclass{30H10, 30H35, 46A55, 47B32, 47B35.}
\thanks{This research was supported in part by grant PID2024-160033NB-I00 funded by MICIU/AEIMCIN/AEI/10.13039/501100011033 and by FEDER, UE}

\begin{document}
\begin{abstract}
Given a reproducing kernel Hilbert space $\mathcal H$ and a bounded operator $T$ on $\mathcal H$, we say that $T$ obeys the {\it strong reproducing kernel thesis} (SRKT) if $\|T\|=\sup_k\|Tk\|$, where $k$ ranges over the normalized reproducing kernels in $\mathcal H$. Two problems are posed. One of these is to determine which Hankel operators on the Hardy space $H^2$ obey the SRKT. The other concerns the geometry of the unit ball of $\text{\rm BMO}$ with respect to the Garsia norm.
\end{abstract}

\maketitle

\section{Introduction}

Suppose $\mathcal H$ is a reproducing kernel Hilbert space. To be more explicit, we shall assume that $\mathcal H$ consists of (complex-valued) functions living on some domain $G\subset\C$ and that the point evaluations $f\mapsto f(z)$ are continuous functionals on $\mathcal H$, for all $z\in G$. The {\it reproducing kernel} associated with a point $z\in G$ is, by definition, the (unique) function $K_z\in\mathcal H$ such that $\langle f,K_z\rangle=f(z)$ for every $f\in\mathcal H$. This last formula implies, in particular, that 
$$\|K_z\|^2_\mathcal H=\langle K_z,K_z\rangle=K_z(z),$$
and we may consider the {\it normalized reproducing kernels} 
\begin{equation}\label{eqn:normrepker}
k_z:=K_z/\sqrt{K_z(z)},\qquad z\in G,
\end{equation}
satisfying $\|k_z\|_\mathcal H=1$. 

\par Suppose further that $\mathcal T$ is a family of linear operators going from $\mathcal H$ to some other Hilbert (or Banach) space, say $\mathcal E$. One says that the operators in $\mathcal T$ obey the {\it reproducing kernel thesis} (RKT, for short) if the condition 
$$\sup\left\{\|Tk_z\|_\mathcal E:\,z\in G\right\}<\infty$$
implies, for $T\in\mathcal T$, the boundedness of $T$ and if this supremum is comparable to $\|T\|_{\mathcal H\to\mathcal E}$. The RKT---when available---thus provides a useful boundedness criterion, and its validity has been established in a number of situations. For instance, the RKT is known to hold for the embedding operators associated with Carleson measures; see, e.g., \cite[Lecture VII]{N}. Another example is provided by Hankel operators (see below).

\par Now, given a bounded linear operator $T:\mathcal H\to\mathcal E$, we say that $T$ obeys the {\it strong reproducing kernel thesis} (SRKT) if 
$$\|T\|_{\mathcal H\to\mathcal E}=\sup\left\{\|Tk_z\|_\mathcal E:\,z\in G\right\}.$$
We shall be interested in this phenomenon for Hankel operators on Hardy spaces of the unit disk, so we proceed to describe the setup.

\par We write $\D$ for the disk $\{z\in\C:|z|<1\}$, $\T$ for its boundary, and $m$ for the normalized arc length measure on $\T$. For $1\le p\le\infty$, the Lebesgue space $L^p:=L^p(\T,m)$ (of complex-valued functions) is introduced in the usual way and equipped with the standard norm $\|\cdot\|_p$. We now define the {\it Hardy space} $H^p$ as the subspace of $L^p$ (carrying the same norm) formed by the functions $f$ whose {\it Poisson integral} 
\begin{equation}\label{eqn:poissint}
\mathcal Pf(z):=\int_\T f(\ze)\,\f{1-|z|^2}{|\ze-z|^2}\,dm(\ze),\qquad z\in\D,
\end{equation}
is analytic on $\D$. Given a function $f\in H^p$, we routinely identify it with its Poisson extension (2) into the disk; the extended function is then still denoted by $f$. Also, $H^p_0$ will stand for the set of $H^p$ functions vanishing at the origin. Our operators will be acting from the Hilbert space $H^2$ to its antianalytic counterpart
$$H^2_-:=L^2\ominus H^2\,(=\ov{H^2_0}),$$
and we write $P_-$ for the orthogonal projection from $L^2$ onto $H^2_-$. 

\par Finally, we recall that $\bmo:=\bmo(\T)$, the space of functions of {\it bounded mean oscillation} on $\T$, can be defined as the set of all functions $h\in L^2$ for which the quantity 
$$\|h\|_G:=\sup_{z\in\D}\left\{\mathcal P(|h|^2)(z)-|\mathcal Ph(z)|^2\right\}^{1/2}$$
(known as the {\it Garsia norm} of $h$) is finite. The subspaces
$$\bmoa:=\bmo\cap H^2\quad\text{\rm and}\quad\bmoa_0:=\bmo\cap H^2_0$$
will also show up in what follows. The reader is referred to \cite[Chapters II and VI]{G} for alternative definitions and basic properties of Hardy spaces and $\bmo$. 

\par Given a function $\ph\in L^2$, the {\it Hankel operator} $H_\ph$ with {\it symbol} $\ph$ is defined (initially on $H^\infty$) by 
$$H_\ph g:=P_-(\ph g),\qquad g\in H^\infty.$$
The classical Nehari theorem (see, e.g., \cite[p.\,181]{N}) states that $H_\ph$ extends to a bounded operator from $H^2$ to $H^2_-$ if and only if it admits a bounded symbol, meaning that $H_\ph=H_\eta$ for some $\eta\in L^\infty$, in which case 
$$\|H_{\ph}\|_{H^2\to H^2_-}=\inf_\eta\|\eta\|_\infty,$$
the infimum being taken over all such $\eta$'s. 

\par It is clear that $H_\ph$ depends only on $P_-\ph(=:\ph_-)$, and moreover, $H_\ph=H_{\ph_-}$. Consequently, every Hankel operator has an {\it antianalytic symbol} (i.e., a symbol from $H^2_-$), which is in fact uniquely determined by the operator. Thus, it seems natural to start with a function $\psi\in H^2_0$ and consider the associated Hankel operator $H_{\ov\psi}$. Because the projection $P_-$ maps $L^\infty$ onto $\ov{\text{\rm BMOA}_0}$, Nehari's theorem implies that $H_{\ov\psi}$ extends to a bounded operator from $H^2$ to $H^2_-$ if and only if $\psi\in\text{\rm BMOA}_0$; it also tells us that
\begin{equation}\label{eqn:hanknorm}
\|H_{\ov\psi}\|_{H^2\to H^2_-}=\|\psi\|_*,
\end{equation}
where
$$\|\psi\|_*:=\inf\left\{\|\eta\|_\infty:\,\eta\in L^\infty,\,\,P_-\eta=\ov\psi\right\}.$$
Next, we look at the action of $H_{\ov\psi}$ on the normalized reproducing kernels $k_z$ in $H^2$. These are related to the Cauchy kernels $K_z(\ze):=(1-\ov z\ze)^{-1}$ as in (1), so that
$$k_z(\ze):=\f{(1-|z|^2)^{1/2}}{1-\ov z\ze},\qquad z\in\D.$$
Furthermore, it is well known that
\begin{equation}\label{eqn:elegiden}
\|H_{\ov\psi}k_z\|_2^2=\mathcal P(|\psi|^2)(z)-|\psi(z)|^2,\qquad z\in\D.
\end{equation}
This elegant identity was discovered by Bonsall in \cite{B}. From (4) it follows at once that
\begin{equation}\label{eqn:hankrkgars}
\sup\left\{\|H_{\ov\psi}k_z\|_2:\,z\in\D\right\}=\|\psi\|_G,
\end{equation}
a fact that was further used in \cite{B} to establish the RKT for Hankel operators. An alternative proof of the RKT in this setting can be found in \cite{T} (in particular, identity (4) appears there as Lemma 2.1).

\par Because the left-hand side of (5) is obviously dominated by that of (3), we have 
\begin{equation}\label{eqn:psitwonorms}
\|\psi\|_G\le\|\psi\|_*.
\end{equation}
One of our problems concerns the SRKT for Hankel operators and reduces to determining the cases of equality in (6). The other, which also involves the Garsia norm, deals with the geometry of the unit ball of $\bmo$ with respect to that norm. Both problems are posed and discussed in the next section.

\section{Problems and discussion}

\medskip\noindent\textbf{Problem 1.} Characterize the Hankel operators (from $H^2$ to $H^2_-$) that obey the SRKT.

\medskip Writing the Hankel operator in question as $H_{\ov\psi}$, with $\psi\in\text{\rm BMOA}_0$, we thus seek to determine whether 
$$\|H_{\ov\psi}\|_{H^2\to H^2_-}=\sup\left\{\|H_{\ov\psi}k_z\|_2:\,z\in\D\right\}.$$
In view of (3) and (5), this boils down to describing the functions $\psi\in\text{\rm BMOA}_0$ for which
\begin{equation}\label{eqn:asteqgars}
\|\psi\|_*=\|\psi\|_G.
\end{equation}

\par Of special interest is the case where $\psi$ is bounded, so that $\psi\in H^\infty_0$. The functions $\eta\in L^\infty$ with $P_-\eta=\ov\psi$ are then precisely those of the form $\ov\psi-h$, where $h$ ranges over $H^\infty$. Consequently, 
\begin{equation}\label{eqn:puzo}
\|\psi\|_*=\inf\left\{\|\ov\psi-h\|_\infty:\,h\in H^\infty\right\},
\end{equation}
and this quantity obviously does not exceed $\|\psi\|_\infty$. We have therefore 
\begin{equation}\label{eqn:bryukho}
\|\psi\|_G\le\|\psi\|_*\le\|\psi\|_\infty.
\end{equation}
In particular, if $\|\psi\|_G$ happens to coincide with $\|\psi\|_\infty$, then (7) holds true and so $H_{\ov\psi}$ obeys the SRKT. 

\par It should be noted, in this connection, that a detailed discussion of the functions $f\in H^\infty$ satisfying $\|f\|_G=\|f\|_\infty$ can be found in \cite{DJFA}; such functions were termed {\it $G$-extremal} there. For instance, every nonconstant inner function is $G$-extremal, and many more examples can be furnished. In fact, for every $g\in H^\infty$ one can find a Blaschke product $B$ such that $Bg$ is $G$-extremal (see \cite[Theorem 2.3]{DJFA} for a more general result). Now, if $f$ ($\in H^\infty$) is $G$-extremal and if $\psi$ ($\in H^\infty_0$) is defined by $\psi=zf$, then $\psi$ is again $G$-extremal, i.e.,
\begin{equation}\label{eqn:psigekstr}
\|\psi\|_G=\|\psi\|_\infty.
\end{equation}
(One may use the elementary inequality 
\begin{equation}\label{eqn:eleminequal}
\|f\|_G\le\|zf\|_G
\end{equation}
to check this.) As we have just seen, (10) implies (7) and ensures the validity of the SRKT for $H_{\ov\psi}$. 

\par Letting again $\psi=zf$, with $f\in H^\infty$, we go on to rewrite the infimum from (8) as 
$$\inf\left\{\|\ov f-h\|_\infty:\,h\in H^\infty_0\right\}=:I(f).$$
We then plug the resulting expression for $\|\psi\|_*$ into (9) and combine it with (11) to get 
\begin{equation}\label{eqn:fiff}
\|f\|_G\le I(f)\le\|f\|_\infty.
\end{equation}
The functions $f\in H^\infty$ satisfying $I(f)=\|f\|_\infty$ were recently brought into view and studied, under the name of {\it saturated} functions, by Brevig and Seip in \cite{BS}. From (12) it is clear that every $G$-extremal function is saturated. Under certain assumptions---and for less obvious reasons---the converse turns out to hold as well. In particular, the two properties are equivalent within the class of nonconstant functions whose outer factor is in the disk algebra $H^\infty\cap C(\T)$; this can be verified by juxtaposing \cite[Theorem 3]{BS} with \cite[Corollary 2.6]{DJFA}.

\par On the other hand, we wonder if there are unbounded functions $\psi\in\text{\rm BMOA}_0$ satisfying (7). Examples, if existent, would be welcome.

\par To conclude our discussion of Problem 1, we remark that there definitely exist Hankel operators which do {\it not} obey the SRKT, since there exist functions $\psi\in H^\infty_0$ with $\|\psi\|_G<\|\psi\|_*$. In fact, a direct calculation shows that this last inequality holds when 
$\psi(z)=cz+z^2$, with a suitable $c>0$. (I owe this observation to Artur Nicolau.) See the next section for a larger class of examples.

\par By way of preparation for our second problem, we now recall that, given a convex set $S$ in a vector space $X$, a point $x$ of $S$ is said to be its {\it extreme point} if the only vector $y\in X$ satisfying $x+y\in S$ and $x-y\in S$ is $y=0$. The convex sets to be dealt with below are 
$$\mathcal B_G:=\left\{f\in\bmo:\,\|f\|_G\le 1\right\}$$
and $\mathcal B_{G,A}:=\mathcal B_G\cap\bmoa$, that is, the unit balls of $\bmo$ and $\bmoa$ with respect to the Garsia norm.

\medskip\noindent\textbf{Problem 2.} Describe the extreme points of $\mathcal B_G$ and of $\mathcal B_{G,A}$.

\medskip Here $\bmo$, as well as $\bmoa$, is actually viewed as a quotient space modulo $\C$; its elements are thus cosets of the form $f+\C$ rather than individual functions $f$. (Since the Garsia norm of a constant function is $0$, this convention is needed for $\|\cdot\|_G$ to become a true norm.) Alternatively, to do without cosets, we may replace $\bmo$ and $\bmoa$ by 
$$\bmo_0:=\left\{f\in\bmo:\,\int_\T f\,dm=0\right\}$$
and $\bmoa_0$, respectively.

\par In \cite{DJFA}, we gave a sufficient condition for a function to be an extreme point of $\mathcal B_G$ (or of $\mathcal B_{G,A}$). Let us say, as we did there, that a function $f\in\bmo$ is 
{\it norm-attaining} if there exists a point $z_0\in\D$ such that 
$$\|f\|_G^2=\mathcal P(|f|^2)(z_0)-|\mathcal Pf(z_0)|^2.$$
Now \cite[Theorem 2.7]{DJFA} states that if $f$ is a norm-attaining function with $\|f\|_G=1$, then $f$---or, strictly speaking, the corresponding coset---is an extreme point of $\mathcal B_G$. The question of whether the converse is true remains open. 

\par For certain natural norms on $\bmo$ other than $\|\cdot\|_G$, the extreme points of the unit ball that arises were studied in \cite{AS} and in \cite{K}. The last-mentioned paper, specifically Theorem 4.1 in \cite{K}, actually characterizes the extreme points of the unit ball in $L^\infty/H^\infty$, a space which is isometrically isomorphic to $\bmoa_0$ endowed with the norm $\|\cdot\|_*$. A similar method was recently employed in \cite{DAAMP} to treat some other $H^p$-related quotient spaces. 

\par In conclusion, we suggest considering an analog of Problem 2 that asks for a description of the {\it strongly extreme} points of $\mathcal B_G$ and/or $\mathcal B_{G,A}$. (By definition, given a Banach space $X$ and a point $x$ in its closed unit ball, $x$ is strongly extreme for the ball if it has the following property: for every $\eps>0$ there is a $\de>0$ such that the inequalities $\|x\pm y\|<1+\de$ imply, for $y\in X$, that $\|y\|<\eps$.) We refer to \cite{CT} and \cite{DAFM} for a discussion of strongly extreme points in some classical function spaces. We also cite \cite{DIEOT}, where further open questions about extreme points can be found. 

\section{Addendum}

After the writing had been essentially completed, a further piece of information on Problem 1 was brought to my attention by Ole Brevig. The following facts, which I learned from him, are included here with his kind permission. 

\par The norm-attaining functions $\psi\in\text{\rm BMOA}_0$ for which the Hankel operator $H_{\ov\psi}$ obeys the SRKT are all of the form
\begin{equation}\label{eqn:innbrev}
\psi=c\left(I-I(0)\right),
\end{equation}
where $c\in\C$ and $I$ is an inner function that has a zero in $\D$. More precisely, if 
$$\|\psi\|^2_G=\mathcal P(|\psi|^2)(z_0)-|\psi(z_0)|^2$$
with a specific $z_0\in\D$, then (7) is equivalent to (13), where the inner function $I$ satisfies $I(z_0)=0$. (For $z_0=0$, this is ensured by \cite[Theorem 1]{Bre}. The general case can be reduced to this by means of a suitable M\"obius transformation that sends $z_0$ to $0$.) 

\par In particular, the SRKT problem is thereby solved for compact Hankel operators $H_{\ov\psi}$. Indeed, the functions $\psi$ that arise in this case are in $\vmo$ and hence norm-attaining. Yet another nice corollary is that, among the polynomials in $H^2_0$, the only functions $\psi$ satisfying (7) are the monomials.

\par Thus, in both problems (i.e., in Problems 1 and 2 above), it is the case of a non-norm-attaining function that remains to be settled.

\section*{Acknowledgements}

I thank Ole Brevig for his illuminating remarks on the first draft of this note. I am also indebted to Artur Nicolau for a helpful discussion, and to Diana Serrano Rodr\'iguez and Daniel N\'u\~nez Alarc\'on for stimulating correspondence.

\section*{Competing Interests}

The author has no competing interests to declare that are relevant to the content of this article.

\section*{Data availability}

This work has no associated data.


\begin{thebibliography}{1} 

\bibitem{AS} Axler, S., Shields, A.: Extreme points in VMO and BMO. Indiana Univ. Math. J. \textbf{31}(1), 1--6 (1982)

\bibitem{B} Bonsall, F.F.: Boundedness of Hankel matrices. J. London Math. Soc. \textbf{29}(2), 289--300 (1984)

\bibitem{Bre} Brevig, O.F.: Minimal norm Hankel operators. Proc. Amer. Math. Soc. \textbf{150}(10), 4383--4391 (2022)

\bibitem{BS} Brevig, O.F., Seip, K.: Maximal norm Hankel operators. J. Math. Anal. Appl. \textbf{529}(2), 127221 (2024)

\bibitem{CT} Cima, J.A., Thomson, J.: On strong extreme points in $H^p$. Duke Math. J. \textbf{40}, 529--532 (1973)

\bibitem{DAFM} Dyakonov, K.M.: Inner functions as strongly extreme points: stability properties. Ann. Fenn. Math. \textbf{48}(2), 681--690 (2023)

\bibitem{DIEOT} Dyakonov, K.M.: Questions about extreme points. Integral Equ. Oper. Theory \textbf{95}(2), 14 (2023)

\bibitem{DJFA} Dyakonov, K.M.: Extremal problems in BMO and VMO involving the Garsia norm. J. Funct. Anal. \textbf{288}(7), 110833 (2025)

\bibitem{DAAMP} Dyakonov, K.M.: Extreme points in quotients of Hardy spaces. Anal. Math. Phys. \textbf{16}(5), 110 (2026)

\bibitem{G} Garnett, J.B.: Bounded Analytic Functions. Revised first edition. Springer, New York (2007)

\bibitem{K} Koosis, P.: Weighted quadratic means of Hilbert transforms. Duke Math. J. \textbf{38}, 609--634 (1971)

\bibitem{N} Nikol'ski\u\i, N.K.: Treatise on the Shift Operator, Grundlehren der mathematischen Wissenschaften, 273. Springer-Verlag, Berlin (1986)

\bibitem{T} Treil, S.: A remark on the reproducing kernel thesis for Hankel operators. St. Petersburg Math. J. \textbf{26}(3), 479--485 (2015)

\end{thebibliography}
\end{document}